\documentclass[12pt]{amsart}
\usepackage{graphicx} 

\usepackage{amsmath,amssymb}

\usepackage{caption}

\title[Tropical Limit of Hyperbolic Amoebas of Surfaces]{Tropical limit of hyperbolic amoebas\\ of complex analytic surfaces}
\author{Peter Petrov and Mikhail Shkolnikov}
\date{}

\newcommand{\SL}{\operatorname{SL}_2(\mathbb{C})}
\newcommand{\SU}{\operatorname{SU}(2)}
\newcommand{\Hy}{\mathbb{H}^3}

\newtheorem*{theorem*}{Theorem}
\newtheorem*{lemma*}{Lemma}

\begin{document}
\begin{abstract}
    In this letter, we establish a general fact about the convergence of images of families of closed analytic surfaces in the special linear group $\SL$ under the quotient by its maximal compact subgroup $\SU$ subject to a contracting scaling sequence. 
\end{abstract}
\maketitle

\section{The result}
Consider $\SL,$ the group of complex two-by-two matrices with unit determinant, and its maximal compact subgroup $\SU$ consisting of all two-by-two unimodular unitary matrices. The homogeneous quotient manifold $\SL\slash \SU$ is contractible and is endowed with a natural complete Riemannian metric of constant sectional curvature $-1,$ therefore modeling the hyperbolic three-space $\Hy$ (see \cite{bryant1987surfaces,thurston1997three}). Let $\varkappa\colon \SL\rightarrow\Hy$ denote the corresponding quotient map and let $O\in\Hy$ be the image of the identity matrix under $\varkappa$.

For positive $s,$ let $R_s\colon\Hy\rightarrow\Hy$ denote a diffeomorphism performing a rescaling centered at $O$ with factor $s^{-1},$ i.e. $R_s(O)=O$ and $R_s(P)=\gamma(s^{-1}t)$ for any other $P\in\Hy$ having distance $t>0$ from $O,$ where $\gamma\colon\mathbb{R}\rightarrow\Hy$ is a natural parametrization of a unique geodesic passing through $O$ and $P$ such that $\gamma(0)=O$ and $\gamma(t)=P.$

A complex analytic surface $V$ in $\SL$ is a subset such that for every point $Z\in V$ there exist an open neighborhood $U\subset\SL$ of $Z$ and a holomorphic function $F\colon U\rightarrow\mathbb{C}$ satisfying $V\cap U=F^{-1}(0).$ Our result (announced previously in \cite{shkolnikov2024introduction} as Claim 1) is the following.

\begin{theorem*}
Let $\{s_n\}_{n=1}^{\infty}$ be a sequence of positive real numbers such that $s_n\rightarrow+\infty,$ and let $\{V_n\}_{n=1}^{\infty}$ be a sequence of closed analytic complex surfaces inside $\SL$. Then, either the sequence $R_{s_n}\circ\varkappa(V_n)\subset\Hy$ goes to infinity, or it has Hausdorff convergent subsequences with the limit of each being a complement to an open ball centered in $O.$
\end{theorem*}

Section \ref{sec_proof} is dedicated to the formal proof of this theorem, which involves standard arguments on complex intersections and Hausdorff metric, as well as an explicit analytic estimate; and  Section \ref{sec_arg} gives a pure geometric explanation for why the statement should be trusted.

\section{The context}
In this section, we justify the title of the article by providing a slightly broader perspective on the problem at hand and clarifying the terms.

First of all, the notion ``amoeba'' stems from an extremely influential book \cite{gelfand1994discriminants}. In that context, it refers to an image of an algebraic subvariety of the complex algebraic torus $(\mathbb{C}\backslash\{0\})^n$ under the coordinate-wise logarithm of the absolute value. In principle, one could take closed analytic subvarieties instead of algebraic ones -- some of the basic facts, such as convexity of connected components of the complement to an amoeba of a hypersurface, will hold. However, the literature on this subject is very scarce (the only publication we are aware of is \cite{madani2015analytic}). 

Despite this, we do believe that this direction is of paramount importance. A particular amoeba of fundamental interest we would like to mention is one of a curve parametrized by $z\rightarrow(\exp(i(z-\frac{1}{2})^2),\zeta(z)).$ To figure out the number of corresponding connected components is (more or less) equivalent to the understanding the Riemann hypothesis.

The map $\varkappa$ defined in the previous section, is the analogue of the logarithmic projection -- both maps are equivalents of taking the quotient by the maximal compact subgroup. Therefore, for a subvariety of $\SL,$ its image in $\Hy$ under $\varkappa$ is called its ``hyperbolic amoeba''. In \cite{mikhalkin2022non}, it was shown that a complement to a hyperbolic amoeba of an algebraic surface is $\SL$ is convex -- in particular, it consists of at most one connected component. Notably, the argument presented there doesn't extend to analytic surfaces since global intersection theory is applied (and the number of isolated intersections of subvarieties of complementary dimensions, in the analytic world, can be infinite).

In the same paper, tropical limits of hyperbolic amoebas of algebraic curves were studied. Here, the tropical limit is exactly what was considered in the previous section -- i.e. we apply some scaling sequence to the amoebas and look at convergent subsequences -- in complete analogy with classical geometric approach to tropicalization (see for instance \cite{mikhalkin2018quantum}, where scaled sequences are used to define tropical limit). The main result in this direction comprised of restricting the set of possible limits claiming, in particular, that they consist of a finite number of spheres centered at $O$ and geodesic segments (whose extensions pass through $O$). Similarly to tropical geometry in $\mathbb{R}^3$, we don't know which of these sets in $\Hy$ can be realized as tropical limits of complex curves.

On the other hand, as we will show below, the tropical limits of hyperbolic amoebas of closed surfaces (even analytic ones!) are easy to describe -- they are complements to open balls centered in $O$. Moreover, we provide examples of realizing every such set as tropical limit.
 
\section{A plausible argument}\label{sec_arg}
A version of the geometric idea reviewed in this section was first presented by the second author at the conference ``Toric and Tropical Techniques in Symplectic Field Theory'' held in June 2024 at the Mittag-Leffler Institute (the recording of the talk is available on the official website of the event -- what was referred to as a ``clean proof'', although in the case of a very similar group $\operatorname{PSL}_2(\mathbb{C}),$ was explained within a couple of minutes starting at 38:15). 

A very useful principle in complex geometry is that a transverse intersection of complex subvarieties of complementary dimension has always multiplicity one, which is due to the presence of canonical orientation coming from holomorphic coordinates (in contrast with real geometric topology where the sign can be negative if orientations do not match). Moreover, even if an isolated intersection is not transverse, the intersection multiplicity is still a positive integer, which can be computed either geometrically (by a small perturbation and reducing to transverse intersections), or topologically (via linking numbers on a small sphere around the intersection point of cycles cut by the subvarieties) or algebraically (as a dimension of the quotient of the ring of power series at the intersection point by the ideal generated by local equations of the subvarieties). 

As a consequence, such an intersection, if we move continuously one of the closed subvarieties, cannot cancel with another intersection, and thus, unless it goes to ``infinity'', it cannot disappear. This principle was erroneously applied to hyperbolic amoebas of surfaces and to hyperbolic amoebas of specific curves \emph{after} the tropical limit.

In more detail, what can be extracted from \cite{mikhalkin2022non} is that lines in $\SL$ may tropicalize to two types of sets. Limit of type A is a set of the form $\{\gamma(t):t\geq r\},$ where $r\geq 0$ and $\gamma$ is a natural parametrization of a geodesic such that $\gamma(0)=O$. Limit of type B is the union of a set of type A, i.e. $\{\gamma(t):t\geq r\},$ and a sphere of radius $r$ centered in $O$.

Type A sets may be used to demonstrate that a complement to a tropical limit (denote it $\Lambda$) of hyperbolic amoebas of a family of surfaces is a star-shaped domain centered at $O.$ Indeed, if it is not, we would be able to find a one parametric family of type A sets with the same $\gamma$ and varying $r$ for which there is a common interval with $\Lambda,$ which shrinks as $r$ approaches $r_0.$ In a similar way, type B sets may be used to show that the function assigning to a direction the distance from $O$ to $\Lambda$ (in this direction) is constant. In both cases, the ``lost'' intersection didn't run to infinity, as it was confined to a bounded region. This completes our original reasoning.

\section{An actual proof}\label{sec_proof}
What was sketched above has some flaws. First, we need to know that every set of type A and B is realizable as a limit of lines, and for every scaling sequence. This, in fact, is easy to verify by a direct computation. Second, we cannot apply the non-cancellation-of-intersection argument after the limit, since there is no single geometric object corresponding to the limiting shape. Thus, we would need to step back, and try to use the rescaled amoebas near the limit. This, in its turn, would require a much finer control of geometry. In addition, we are varying the limits themselves, which assumes that we would need some strong version of uniform convergence. 

What will be described next is a complete and rigorous proof of our Theorem arisen out of attempts to make the above reasoning work. It doesn't accurately follow the former steps, but it still uses hyperbolic amoebas of lines (although without relying on the results on their tropical limit) and local intersection theory.

One of the easy statements we exploit is that a hyperbolic amoeba of a line in $\SL$ is a horosphere in $\Hy$ and every horosphere is realized in this way (the readers might want to either think about this easy statement themselves or consult with \cite{mikhalkin2022non}). Recall that a horosphere is a closed submanifold of $\Hy$ such that it is orthogonal to all geodesics sharing the same point at infinity. All horospheres with the same point at infinity define a trivial fibration of $\Hy,$ the base of which may be identified with a particular geodesic with this point at infinity by assigning to a horosphere the unique intersection point with the geodesic. In particular, a horosphere cuts $\Hy$ into two connected components one of which is convex and the other is not.

The key ingredient in the proof of the Theorem is the following technical fact about two geometric balls and a rescaled horosphere. 

\begin{figure}[h!]
    \centering
    \includegraphics[angle=90,origin=c,width=0.65\linewidth]{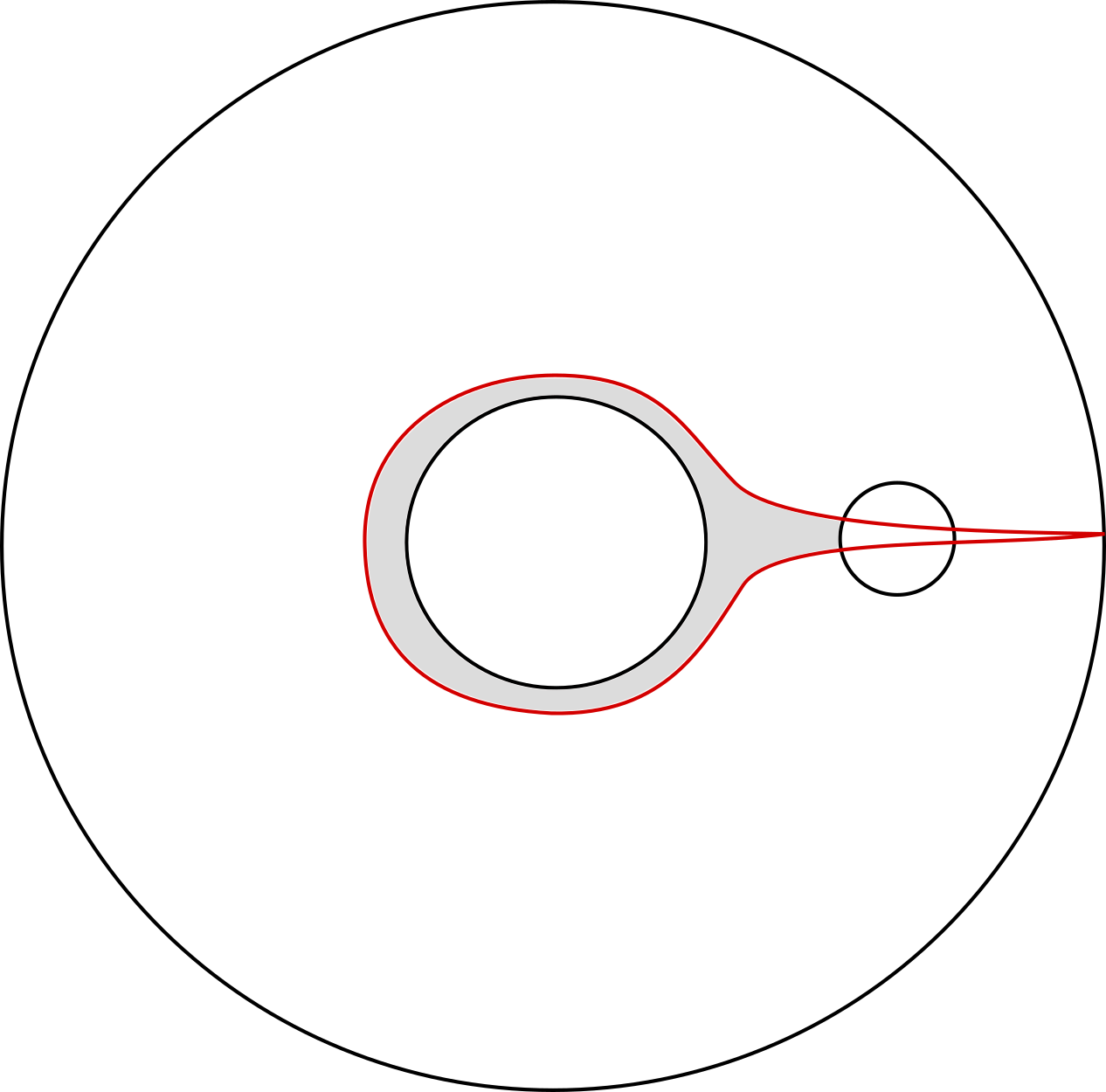}
    \caption{An illustration of the Lemma. A part of rescaled horoshpere together with a part of the boundary of the upper  ball to surround the ball in the middle.}
    \label{fig_horo1}
\end{figure}

\begin{lemma*}
Let $P$ be a point in $\Hy$ such that the distance from $Q$ to $O$ is $d>0.$ In addition let $\varepsilon$ and $\rho$ be positive numbers such that $\varepsilon+\rho>d.$ Then, there exist $\sigma>0$ such that for all $s>\sigma$ there exist a horosphere $H$ such that $R_s(H)\backslash\partial B_{\epsilon}(Q)$ and $\partial B_{\epsilon}(Q)\backslash R_s(H)$ have non-empty connected components $C_1$ and $C_2$ respectively satisfying the property that the closure of their union is a topological sphere and $B_{\rho}(O)$ is contained within the bounded region $D$ surrounded by this sphere, and, in addition the interior of $D$ contains no points of $B_{\epsilon}(Q)$. 
\end{lemma*}

See Figure \ref{fig_horo1} for a schematic illustration of the desired configuration. Our proof of the lemma requires a relatively straightforward  estimate that will be done at the end of this section. 

We will proceed by proving first the main result assuming the Lemma. 

As a general remark, observe that either the sequence $R_{s_n}\circ\varkappa(V_n)$ ``goes to infinity'', i.e. for every open ball $B_r(O)$ its intersection with a member of the sequence is eventually empty, or there exist such $r>0$ that $B_r(O)\cap (R_{s_n}\circ\varkappa(V_n))$ is not empty for infinitely many $n.$ In the former case we may abuse slightly the definition and say that the sequence converges to an empty set (and that the complement of the limit is a ball of infinite radius). 

In the latter case, $\overline{B_r(O)}\cap (R_{s_n}\circ\varkappa(V_n))$ has a convergent subsequence (with non-empty limit) indexed by some $\{n_k\},$ due to the standard fact that the space of non-empty closed sets on a metric compact endowed with Hausdorff metric is a compact itself. From the argument below it will follow that $R_{s_{n_k}}\circ\varkappa(V_{n_k})$ also converges. In fact, we don't really need to rely on the existence of the convergent subsequence, since our argument will imply that. 

\begin{proof}[Proof of the Theorem]

Assume that the sequence $R_{s_n}\circ\varkappa(V_n)$ doesn't run to infinity, in particular there is a point $Q,$ whose neighborhood is visited infinitely often (by a subsequence indexed by $\{n_k\}$). Let $\alpha\geq 0$ be the distance from $Q$ from $O.$ We will show that every $P\in\Hy$ with distance $\beta>\alpha$ from $O$ is a limiting point of $R_{s_{n_k}}\circ\varkappa(V_{n_k})$ (which in its turn implies the Theorem by choosing $Q$ to be the closest limit point to $O$). 
The contrary of this reads as that there exist $\delta>0$ such that $B_\delta(Q)\cap (R_{s_{n_k}}\circ\varkappa(V_{n_k}))$ is empty infinitely often.

Applying the Lemma to a (possibly small) ball $B_1$ around $Q$ and a (possibly big) ball $B_2$ around $O,$ such that the latter contains $P$ in its interior, there exist $s_{n_K}>\sigma$ such that there is a horosphere $H$ with rescaling $R_{s_{n_K}}(H)$ representing the configuration with respect to the mentioned balls such as described in the Lemma. Now, we consider a one-parametric family of contracting horospheres $\{H_t\}_{t\geq 1}$ with the same point at infinity such that $H_1=H.$ 

There exist such $t_0>1$ for which $B_2\cap(R_{s_{n_K}}\circ\varkappa(V_{n_K}))$ intersects $R_{s_{n_K}}(H_{t_0})$ at some point $X$ (which may be thought as being near $P$). At this moment, lift the horosphere $H_{t_0}$ to a line $l_{t_0}$ intersecting $V_{n_K}$ at $Z$ such that $R_{s_{n_K}}\circ\varkappa(Z)=X.$ This is a precise intersection point we are going to loose if we deform the line $l_{t}$ in parallel with its image $H_t$ as $t\geq t_0$ increases, since there are no points of $V_{n_K}$ in $\varkappa^{-1}(B_1).$
\end{proof}

\begin{figure}
    \centering
    \includegraphics[width=0.65\linewidth]{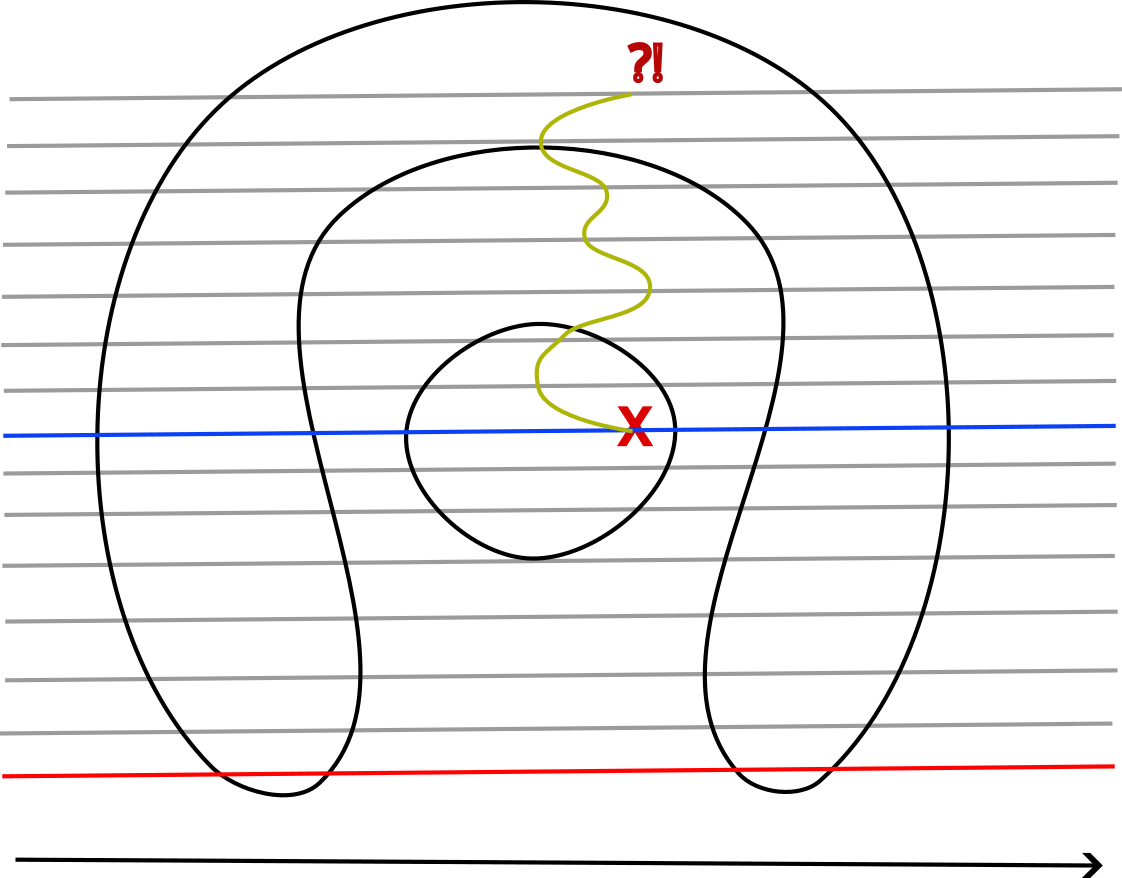}
    \caption{An illustration of the trick used in the proof of Theorem shown in the upper-half space model of $\Hy$ after an inverse of the rescaling $R_{s_{n_K}}$. The region in the middle is where there are some points of a hyperbolic amoeba $\varkappa(V_{n_K}),$ one is marked by X, and the horse-shoe shape around is the region where there are none. Parallel lines represent horospheres moving upwards, and the winding path (stating at X and ending at ?!) corresponds to the hypothetical image of disappearing intersection.}
    \label{fig:enter-label}
\end{figure}

\begin{proof}[Proof of the Lemma]
    Consider the geodesic ray going from $O$ to $Q.$ The horospheres $H$ we are going to construct will be having the end of this ray at infinity as their own limiting point. Note that the whole picture in this presentation becomes invariant with respect to the rotations around the geodesic ray. Therefore, we may prove an an analogous statement for the hyperbolic plane instead of hyperbolic space. In addition, we reverse the perspective: instead of rescaling $H$ we will be doing opposite rescaling of the balls. 

    We use two models of $\mathbb{H}^2.$ One is that of symmetric real positive-definite unimodular matrices -- the distance to $O$ (i.e. the identity matrix) is computed as the logarithm of the absolute value of an eigenvalue, in particular taking real powers of a matrix corresponds to rescaling with center at $O$. Another model is that of a complex upper-half plane, in which a horocycle (i.e. the intersection of a horosphere in $\Hy$ with $\mathbb{H}^2$) is a line parallel to the real axis. The identification between the two models is done via M\"obius transformation applied to $i=\sqrt{-1}$.

    In the matrix model, consider a point at distance $r>0$ from $O,$ i.e. a matrix $M=e^rv_1v_1^\intercal+e^{-r}v_2v_2^\intercal$ where $v_1,v_2$ is an orthonormal eigenvector basis of $\mathbb{R}^2.$ In particular, $$v_1v_1^\intercal=\begin{pmatrix}a&b\\b& c\end{pmatrix} \text{\ and\ } v_2v_2^\intercal=\begin{pmatrix}c&-b\\-b& a\end{pmatrix},$$ for some real $a,b,c$ such that $a,c>0,$ $ac-b^2=0$ and $a^2+c^2=1.$ Applying the M\"obius transformation associated with $M$ to the point $i$ and taking the imaginary part results in $$\operatorname{Im}\frac{(e^ra+e^{-r}c)i+(e^rb-e^{-r}b)}{(e^rb-e^{-r}b)i+(e^rc+e^{-r}a)}=\frac{1+ac}{e^{2r}(c^2+ac)+e^{-2r}(a^2+ac)}.$$
    This formula has several consequences. First, unless $a=c=0$ (which corresponds to a point on the ray we considered at the beginning of this proof), the imaginary part of a point associated with $M^s$ decays asymptotically as $\operatorname{const}\cdot e^{-2rs}$ when $s\rightarrow+\infty,$ with a uniform lower bound being $e^{-2rs}$ itself. 

    Now, consider a point $S$ at the boundary of the ball $B_\varepsilon(Q)$ which doesn't belong to the imaginary ray in the complex picture. This point $S$ has a distance from $O$ strictly greater than $\rho.$ In particular, the asymptotic of the imaginary part of the complex counterpart of $S^s$ is strictly smaller than $e^{-2\rho s}.$ This implies that there exist $\sigma>0,$ such that for all $s>\sigma$ the minimal value (denote it $\mu_s$) of the imaginary part of the $s$-rescaled boundary of $B_\varepsilon(Q)$ is lower than the imaginary part of the $s$-rescaled boundary of $B_\rho(O).$ The desired horocycle is given by the horizontal line $\operatorname{Im}=\mu_s.$ 
\end{proof}

\section{Remarks}
We decided that this short text is not a place to give a detailed account of foundations of tropical geometry and its significance (some basics, which are not used here anyway, may be found in \cite{brugalle2015brief} and, a biased overview of some applications is a major part of the introduction of our previous article \cite{shkolnikov2024introduction}). In the same spirit, we didn't try to inform the reader about the fundamental importance of the groups $\SL$ and $\SU$ in mathematics and physics, as well as of the corresponding quotient space $\Hy$ and its ambient Minkowski space. We only share that to our surprise, the $\SL$ version of the logarithmic projection, the map $\varkappa$ in our notation, was utilized by differential geometers since at least \cite{bryant1987surfaces} (and later \cite{galvez2000flat}), where it was used to describe constant-mean-curvature-$1$ (and, respectively, flat) surfaces in $\Hy$ as images of special types of holomorphic curves in $\SL.$

As was mentioned, the theorem stated in the first section complements the result on tropical limits of hyperbolic amoebas of curves \cite{mikhalkin2022non}, which still is not ultimate since we only know that such limits are described by certain types of floor diagrams but it is unclear which of these diagrams are realized as limits. In the present article the situation is better  -- each complement of an open ball of radius $r$ centered in $O$ is realized by a tropical limit. To see this, consider for instance a family of surfaces $V_n=\{A\in\SL:\operatorname{Tr}(A)=n^r\}$ with the scaling sequence $s_n=\log(n).$ In addition, if instead of $n^r$ we would write $\exp(n),$ the rescaled sequence of amoebas would go to infinity, i.e. we may say that the limit is empty. It is also worth mentioning that for a constant family of closed surfaces, its tropical limit is the whole $\Hy.$

A striking difference of our result with the grand majority of the existing literature on amoebas and tropical geometry is that the presented approach works equally well for algebraic and analytic surfaces. For instance, it was known already from \cite{mikhalkin2022non} that the complement to a hyperbolic amoeba of an algebraic surface is convex (and in particular, connected -- although, it can be unbounded) -- however, the technique used there relies on global intersection theory in the compactification, which is not available for non-algebraic surfaces. In particular, although we don't know at the moment if there is at most one connected component in the complement of a hyperbolic amoeba of a closed non-algebraic surface, in the limit only one component is visible. An analytic way to show that there is still only one component before the limit, would be via extending the notion of Ronkin function to the non-commutative setup and proving that it is locally constant on the complement to an amoeba and geodesically concave on whole $\Hy$.

Interestingly enough, the rescaling diffeomorphism $R_s$ may spoil the convexity of a domain in $\Hy,$ in contrast with homotheties in a Euclidean space. For instance, one of the connected components in the complement to a horosphere is convex, but after a rescaling it may no longer be true as it is seen on Figure \ref{fig_horo1} (we may imagine that the drawing is made in the Klein model where the two notions of convexity, visual and hyperbolic, coincide). Our remark here is that nevertheless, the limit of hyperbolic amoebas of surfaces has a convex complement (a ball), even though before the limit it may be non-convex.

It is instructive to contemplate the meaning of our proof at a less geometric level. Assume that the family of analytic surfaces $\{V_n\}$ in $\SL$ is given by a sequence of equations $f_n=0,$ where each $f_n$ is a global holomorphic function. What we are saying is that if one finds a sequence of solutions $A_n\in \SL$ (i.e. satisfying $f_n(A_n)=0$) such that the norm $|A_n|$ goes to infinity (the norm is the square root of the sum of squares the absolute values of entries), then for any prescribed one-dimensional subspace $L$ of $\mathbb{C}^2$ and a number $\lambda\geq 1$ there exist a sequence $B_n\in\SL$ such that $f_n(B_n)=0,$ $\log|B_n|\sim\lambda\log|A_n|$ and the rescaling $|B_n|^{-1}B_n$ has only rank one matrices with image $L$ as its accumulation points. To derive this statement, we have used $s_n=\log |A_n|$ as the scaling sequence. Alternatively, we may find a sequence of solutions $B_n$ with prescribed kernel of the limit matrices.

A valid question now is whether both the image and the kernel of limit points of $|B_n|^{-1}B_n$ can be prescribed independently. The answer is ``no'', at least in the algebraic case and for arbitrary $r$, which is explained by the structure of $\operatorname{PSL}_2$ phase tropicalization introduced in  \cite{shkolnikov2024introduction} and further developed in \cite{shkolnikov2024PSL2} -- only for some {\it critical} levels $r$ it is possible, and otherwise there is an algebraic relation between the slopes of the image and the kernel lines. In fact, the Theorem proven in this letter served as the motivation to restore the phase, in which case we start seeing a finer structure over the complement to the ball -- over a point $P\in\Hy$ at regular level there is a union of geodesic circles in $\varkappa^{-1}(P)\cong S^3$ and at critical levels one observes two-dimensional fibers with boundary on these circles. This description of a phase diagram, however, is still preliminary, and rather corresponds to a few specific examples and incomplete statements about some families of surfaces. Provided that we have a full control of all phase diagrams, an algebraic (and weaker) version of the present Theorem could be deduced in future. 
\newpage
 Our result is expected to admit generalizations to other reductive complex groups, and be applied in some of the respective proofs, although a description of the corresponding tropical limits will certainly be more complicated, since the most general statement in this direction should extend both our Theorem for $\SL$ (literally the same fact trivially holds for $\operatorname{PSL}_2(\mathbb{C})$ via the two-fold covering) and classical Kapranov's Theorem for tropicalizations of algebraic hypersurfaces in an algebraic torus, producing tropical hypersurfaces -- polyhedral complexes in $\mathbb{R}^n$ with faces of rational slopes -- which are indeed objects of very different geometric flavor than our complement to a ball in $\Hy$. 

\bibliographystyle{plain}
\bibliography{bibliography}
\vspace{18pt}
\address{\begin{center}
     {Institute of Mathematics and Informatics\\\vspace{3pt} at the Bulgarian Academy of Sciences\vspace{3pt}\\ Akad. G. Bonchev St, Bl. 8, 1113 Sofia, Bulgaria}\end{center}}
\vspace{3pt}
\email{\begin{center}
    pk5rov@gmail.com\\ m.shkolnikov@math.bas.bg \end{center}}

\end{document}